\documentclass[reqno,11pt]{amsart}

\DeclareMathOperator*{\esssup}{\mathrm{ess\,sup}}

\usepackage[a4paper,  margin=3.4cm]{geometry}

\usepackage{amsmath,amssymb}
\usepackage{dsfont}         % \mathds{R} np.
\usepackage{comment}
\newtheorem{theorem}{Theorem}[section]

\newtheorem{proposition}[theorem]{Proposition}
\theoremstyle{definition}
\newtheorem{definition}[theorem]{Definition}
\theoremstyle{remark}

\numberwithin{equation}{section}

\title{Dynamics of an infinite age-structured particle system}

\author{Dominika Jasi\'nska}

\address{Instytut Matematyki, Uniwersytet Marii Curie-Sk{\l}odowskiej, 20-031 Lublin, Poland}
\email{jasdominika@wp.pl}

\author{Yuri Kozitsky}

\address{Instytut Matematyki, Uniwersytet Marii Curie-Sk{\l}odowskiej, 20-031 Lublin, Poland}
\email{jkozi@hektor.umcs.lublin.pl}


\begin{document}

\begin{abstract}

The Markov evolution is studied of an infinite age-structured
population of migrants arriving in and departing from a continuous
habitat $X \subseteq\mathds{R}^d$ -- at random and independently of
each other. Each population member is characterized by its age
$a\geq 0$ (time of presence in the population) and location $x\in
X$. The population states are probability measures on the space of
the corresponding marked configurations. The result of the paper is
constructing the evolution $\mu_0  \to \mu_t$ of such states by
solving a standard Fokker-Planck  equation for this models. We also
found a stationary state $\mu$ existing if the emigration rate is
separated away from zero. It is then shown that $\mu_t$ weakly
converges to $\mu$ as $t\to +\infty$.

\end{abstract}

\maketitle

%\tableofcontents

\section{Introduction}

The stochastic dynamics of structured populations attract
considerable attention, see, e.g., quite recent works
\cite{Fima,Koz,Koz1,bib7,Tit}. This, in particular, relates to
age-structured populations studied at both microscopic or mesoscopic
scales, cf. \cite{Fima,Koz1,bib7} and \cite{Wang}, respectively.
Finite populations of this kind are much more mathematically
accessible in contrast to infinite ones where only few results were
obtained. Mostly because the states of such systems -- probability
measures on the corresponding spaces of infinite configurations --
are quite abstract objects, not appropriate for a direct
investigation.  The present work is a continuation of that in
\cite{bibm} where the study was initiated of infinite age-structured
populations based on the use of correlation functions. By employing
such functions one can deal with infinite systems indirectly. The
main advantage of this approach (see \cite{Koz} and the works quoted
therein) is that correlation functions are defined on spaces of
\emph{finite} configurations, that allows one to employ more
powerful tools of solving corresponding evolution equations.

In this work, we introduce and study an individual-based
(microscopic) model of an infinite particle systems --
age-structured populations of migrants. The population dwells in a
spatial habitat, $X\subseteq\mathbb{R}^d$, $d\geq 0$, and each
population member - entity - is characterized by its spatial
location $x\in X$ and age $a\in \mathbb{R}_{+}:=[0,+\infty)$. The
entities arrive (appear) and depart (disappear) at random --
independently of each other. By $\hat{x} = (x,a_x)$ we denote the
corresponding compound trait. Pure states of the population are
collections $\hat{\gamma}$ (called \emph{configurations}) of the
traits of its members. We assume that neither two of these members
can have the same spatial location. Due to this assumption we can
employ here techniques of the theory of marked configuration spaces,
see \cite{bib6}. The set of all such configurations $\hat{\Gamma}$
is endowed with a topology, and hence with the corresponding Borel
$\sigma$-field of measurable subsets, see below.

As mentioned above, in our model the entities arrive and depart
independently of each other at rates (probability density per time)
$b(x)$ and $m(\hat{x})$, respectively. In view of the random
character of the evolution, the population states are probability
measures on $\hat{\Gamma}$, and their Markov evolution is described
by the Kolmogorov equation
\begin{equation}
  \label{KE}
  \frac{d}{dt} F_t = L F_t, \qquad F_{t}|_{t=0} = F_0,
\end{equation}
where $L$ is supposed to be a model-specific linear operator and
$F:\hat{\Gamma}\to \mathbb{R}$ stands for an observable (test
function). Then the weak evolution of states $\mu_0 \to \mu_t$ is
obtained (Theorem \ref{1Tm}) by solving the Fokker-Planck equation
corresponding to (\ref{KE}), see (\ref{FPE}) below. Under some
additional assumption we find also a stationary state of this
evolution and prove the weak convergence of $\mu_t$ to this state as
$t\to +\infty$. In Section 2, we introduce all necessary notions and
facts, whereas in Section 3 we prove Theorem \ref{1Tm}. In
particular, we construct the evolution $\mu_0 \to \mu_t$ in an
explicit form. In subsequent works, we will use this construction to
study the evolution of similar age-structured populations with
interactions.

\section{Preliminaries}
The state of an entity in the population is characterized by its
compound trait $(x,a_x)$, $x\in X\subseteq \mathds{R}^d$ and $a_x\in
\mathbb{R}_{+}$. We use the following notations $\hat{x}=(x,a) \in
\hat{X} = X \times \mathbb{R}_{+}$. For a function $g: \hat{X}\to
\mathbb{R}$, we use interchangeable writings $g(\hat{x})$ and
$g(x,a)$.

\subsection{Marked configurations}
The pure state of the whole population is the collection
$\hat{\gamma}$ of the traits of all its members. The set of all such
(marked) configurations $\hat{\Gamma}$ is endowed with the topology
which we introduce now. Here we mostly follow the approach of
\cite[Sect. 2]{bib6}. First, we define the underlying configurations
space
\begin{equation}
  \label{CF}
\Gamma=\{ \gamma \subset {\mathbb{R}}^d: |\gamma \cap
\varLambda|<\infty \text{ for a compact} \ \varLambda \subset
{\mathbb{R}}^d \}.
\end{equation}
This space is endowed with the vague (weak-hash) topology, cf.
\cite{DV1}, which is the weakest topology that makes continuous all
the maps $\Gamma \ni \gamma \mapsto \langle \gamma, f \rangle$,
$f\in C_{\rm cs}(X)$. Here
\[
\langle \gamma, f \rangle = \sum_{x\in \gamma} f(x),
\]
and $C_{\rm cs}(X)$ stands for the set of all continuous compactly
supported functions $f:X \to \mathbb{R}$. Along with the space
defined in (\ref{CF}) we also use the space of finite configurations
\begin{equation}
  \label{CS1}
  \Gamma_{0}=\bigcup_{n\in {\mathbb{N}}_0} \{ \gamma \subset
{\mathbb{R}}^d : |\gamma|=n \}.
\end{equation}
For a given $\hat{\gamma} \subset X \times \mathds{R}_{+}$, we set
\[
p(\hat{\gamma}) =  \{x\in X: (x,a_x)\in \hat{\gamma}\}.
\]
Then
\begin{equation}
  \label{CF1a}
\hat{\Gamma}:=\{ \hat{\gamma}: p(\hat{\gamma})\in \Gamma\}.
\end{equation}
According to (\ref{CF1a}) neither of two elements of each
configuration $\hat{\gamma}\in \hat{\Gamma}$ can have the same
spatial location. Let now $\mathcal{C}$ denote the set of bounded
continuous functions $g: X \times \mathbb{R}_{+}\to \mathbb{R}$,
each of which is supported on $\varLambda\times \mathbb{R}_{+}$ for
a compact $\varLambda \subset X$. Then the topology of
$\hat{\Gamma}$ is defined as the weakest topology that makes
continuous all the maps
\[
\hat{\Gamma}\ni \hat{\gamma} \mapsto \sum_{ x\in p(\hat{\gamma})} g
(x, a_x) , \qquad g \in \mathcal{C}.
\]
It is known, see \cite[Lemma 2.1]{bib6}, that this topology is
metrizable in such a way that the obtained metric space  is complete
and separable. Let $\mathcal{B}(\hat{\Gamma})$ stand for the
corresponding Borel $\sigma$-field of subsets of $\hat{\Gamma}$. By
$\mathcal{P}(\hat{\Gamma})$ we denote the set of all probability
measures on $(\hat{\Gamma}, \mathcal{P}(\hat{\Gamma}))$. Next we set
\begin{equation*}
  %\label{CS2}
  \hat{\Gamma}_0 = \{ \hat{\gamma}\in  \hat{\Gamma}:
  p(\hat{\gamma})\in \Gamma_0\},
\end{equation*}
where $\Gamma_0$ is defined in (\ref{CS1}). Then $\hat{\Gamma}_0$ is
endowed with the topology induced by the topology of $\hat{\Gamma}$,
and thus with the corresponding Borel $\sigma$-field
$\mathcal{B}(\hat{\Gamma}_0)$. It can be shown that a function
$G:\hat{\Gamma}_0\to \mathbb{R}$ is measurable if and only if there
exists a collection $\{G^{(n)}\}_{n\in \mathbb{N}_0}$ of symmetric
Borel functions $G^{(n)}: \hat{X}^n \to \mathbb{R}$, $\hat{X}:=
X\times \mathbb{R}_{+}$ such that $G^{(0)} = G(\varnothing)$ and
\begin{equation}
  \label{CS3}
  G(\hat{\gamma}) = G^{(n)} (\hat{x}_1, \dots , \hat{x}_n), \quad
  {\rm for} \ \
  \hat{\gamma}=\{\hat{x}_1, \dots , \hat{x}_n\}, \quad n\in \mathbb{N}.
\end{equation}
For such functions, we set
\begin{eqnarray*}
  %\label{CS4}
\int_{\hat{\Gamma}_0} G(\hat{\gamma}) \hat{\lambda} ( d
\hat{\gamma}) = G(\varnothing) + \sum_{n=1}^\infty\frac{1}{n!}
\int_{\hat{X}^n} G^{(n)} (\hat{x}_1, \dots , \hat{x}_n) d \hat{x}_1
\cdots d \hat{x}_n,
\end{eqnarray*}
where $d \hat{x}$ is the Lebesgue measure on $\mathbb{R}^d \times
\mathbb{R}_{+}$. This defines a locally finite measure
$\hat{\lambda}$ on $(\hat{\Gamma}_0,\mathcal{B}(\hat{\Gamma}_0))$ --
the Lebesgue-Poisson measure.  It readily satisfies
\begin{equation}
  \label{Minlos}
\int_{\hat{\Gamma}_0} \sum_{\xi \subset p(\hat{\eta})} G(\hat{\xi},
\hat{\eta}\setminus \hat{\xi}) \hat{\lambda} (d \hat{\eta}) =
\int_{\hat{\Gamma}_0} \int_{\hat{\Gamma}_0} G(\hat{\xi}, \hat{\eta})
\hat{\lambda} (d \hat{\xi}) \hat{\lambda} (d \hat{\eta}),
\end{equation}
that holds for all appropriate functions.

For a measurable $F:\hat{\Gamma}\to \mathbb{R}$ and some $\mu \in
\mathcal{P}(\hat{\Gamma})$, we write
\[
\mu(F) = \int_{\hat{\Gamma}} F(\hat{\gamma}) \mu (d \hat{\gamma}).
\]
A collection, $\mathcal{F}$, of functions $F:\hat{\Gamma}\to
\mathbb{R}$ is called \emph{separating} (measure-defining) if,  for
any two probability measures, $\mu(F) = \nu(F)$ holding for all
$F\in \mathcal{F}$ implies $\mu=\nu$.
\begin{proposition}
  \label{J1pn}\cite[Theorem 1.3.26, page 113]{Mon}
Let $V$ and $\mathcal{F}$ be a complete and separable metric spaces
and a family of functions $F:V \to \mathds{R}$, respectively. Assume
that: (a) each $F\in \mathcal{F}$ is bounded and continuous; (b) for
$F_1, F_2\in \mathcal{F}$, their pointwise product is in
$\mathcal{F}$; (c) for each distinct $v_1,v_2\in V$, there exists
$F\in \mathcal{F}$ such that $F(v_1)\neq F(v_2)$; (d) $\mathcal{F}$
contains $F\equiv 1$. Then $\mathcal{F}$ is separating.
\end{proposition}
By this statement the collection of function
\begin{equation}
  \label{SF}
  F_\theta (\hat{\gamma}) = \prod_{x\in p(\hat{\gamma})}( 1 +
  \theta(x,a_x)),
\end{equation}
is separating, where measurable $\theta: \hat{X} \to (-1,0]$ are
such that each $\theta (x,a) =0$ whenever $x\in \varLambda^c := X
\setminus \varLambda$ for a compact $\varLambda \subset X \subset
\mathds{R}^d$. Note that each $F_\theta$ is measurable and bounded
-- hence $\mu$-integrable for each $\mu\in
\mathcal{P}(\hat{\Gamma})$. Let $q:\hat{X}\to (0,1)$ be a measurable
function. Then each $\theta_q (\hat{x}) :=
q(\hat{x})\theta(\hat{x})$ has the mentioned properties and the
collection of all such $F_{q\theta}$ can be used to determine the
following notion.
\begin{definition}
  \label{0df}
 For a given measurable $q:\hat{X}\to (0,1)$ and $\mu\in \mathcal{P}(\hat{\Gamma})$, the measure $\mu_q$
 defined by the relation $\mu_q(F_\theta) = \mu(F_{\theta_q})$ is
 called an independent $q$-thinning of $\mu$.
\end{definition}
To illustrate this notion, let us take $\mu=\delta_{\hat{\gamma}}$
-- the Dirac measure with atom $\hat{\gamma}$. Then in state $\mu_q$
each $\hat{x}\in \hat{\gamma}$ is retained in $\hat{\gamma}$ with
probability $q(\hat{x})$.

For $\mu_1, \mu_2 \in \mathcal{P}(\hat{\Gamma})$, their convolution
is defined by the relation
\begin{equation}
 \label{Conv}
  (\mu_1 \star \mu_2) (F) = \int_{\hat{\Gamma}^2}
  F(\hat{\gamma}_1\cup \hat{\gamma}_2) \mu_1 (d \hat{\gamma}_1) \mu_2 (d
  \hat{\gamma}_2),
\end{equation}
that ought to hold for all bounded measurable functions
$F:\hat{\Gamma}\to \mathds{R}$. This, in particular, means that
\begin{equation}
  \label{Conv1}
(\mu_1 \star \mu_2) (F_\theta) = \mu_1 (F_\theta) \mu_2(F_\theta).
\end{equation}

\subsection{Tempered configurations}

In this work, we deal with probability measures on $\hat{\Gamma}$
possessing a certain important property. In view of this, we select
a subset of $\hat{\Gamma}$ -- related to this  property -- and thus
`forget' of the remaining configurations. Let $\psi:X \to
\mathds{R}_{+}$ be: (a) continuous, bounded and strictly positive;
(b) integrable, i.e.,
\[
\int_{X} \psi(x) d x <\infty.
\]
One can take $\psi(x) = e^{-|x|}$ as an example of such a function.
Define
\begin{equation}
  \label{J}
  \Psi (\hat{\gamma}) = \sum_{x\in p(\hat{\gamma})} \psi(x), \qquad \hat{\gamma }\in
 \hat{ \Gamma}.
\end{equation}
Note that $\Psi (\hat{\gamma})$ can take infinite values for some
$\hat{\gamma}$. Then the set of tempered configurations is defined
as
\begin{equation}
  \label{J1}
\hat{\Gamma}_* = \{ \hat{\gamma}\in \hat{\Gamma}: \Psi
(\hat{\gamma}) <\infty \}.
\end{equation}
Similarly as in \cite[subsect. 2.3]{Koz2} we equip this set with the
following metric
\begin{equation}
  \label{J2}
  \upsilon_*(\hat{\gamma}, \hat{\gamma}') = \sup_{g} \left|\sum_{x\in p(\hat{\gamma})}
  g(\hat{x}) \psi(x) - \sum_{x\in p(\hat{\gamma}')} g(\hat{x}) \psi(x)
  \right|,
\end{equation}
where the supremum is taken over the subset of the set of bounded
Lipschitz-continuous functions $C^{BL}(\hat{X})$ consisting of those
$g: \hat{X}\to \mathds{R}$ for which
\[
\sup_{\hat{x}\in \hat{X}} |g(\hat{x})| + \sup_{\hat{x}\neq  \hat{y}
\in \hat{X}} \frac{\left| g(\hat{x}) - g(\hat{y})\right|}{|x - y| +
|a_x - a_y|}\leq 1.
\]
It is possible to prove, cf. \cite[Proposition 2.7]{Koz2},  that the
metric space $(\hat{X}, \upsilon_*)$ is complete and separable. Let
$\mathcal{B}(\hat{\Gamma}_*)$ be the corresponding Borel
$\sigma$-field of such subsets of $\hat{\Gamma}_*$. By Kuratowski's
theorem \cite[Theorem 3.9, page 21]{Part} one then proves that
$\hat{\Gamma}_* \in \mathcal{B}(\hat{\Gamma})$ and
$\mathcal{B}(\hat{\Gamma}_*)$ coincides with the Borel
$\sigma$-field related to the topology on $\hat{\Gamma}_*$ induced
by the vague topology of $\hat{\Gamma}$. This allows one to redefine
each $\mu \in \mathcal{P}(\hat{\Gamma})$ with the property
$\mu(\hat{\Gamma}_*)=1$ as a measure on $(\hat{\Gamma}_*,
\mathcal{B}(\hat{\Gamma}_*))$, see \cite[Corollary 2.8]{Koz2} for
further details.
\begin{definition}
  \label{1adf}
By $\mathcal{P}_*$ we denote the set of probability measures $\mu$
on $(\hat{\Gamma},\mathcal{B}(\hat{\Gamma})$  with the property
$\mu(\hat{\Gamma}_*)=1$.
\end{definition}
As just mentioned, each $\mu\in\mathcal{P}_*$ can be redefined as a
probability measure on $(\hat{\Gamma}_*,
\mathcal{B}(\hat{\Gamma}_*))$, which we assume to be done from now
on. Note that $\mu$ belongs to $\mathcal{P}_*$ if and only if
$\mu(\Psi) <\infty$, cf. (\ref{J1}).
\begin{definition}
\label{1df} By $\mathcal{P}_{**}$ we denote the set of probability
measures $\mu$ on $(\hat{\Gamma},\mathcal{B}(\hat{\Gamma}))$ for
each of which $\mu(F_\theta)$ can be written
 in the following form
 \begin{eqnarray}
   \label{SF1}
   \mu(F_\theta) = \int_{\hat{\Gamma}_0} \hat{k}_\mu(\hat{\gamma})\Big( \prod_{x\in
   p(\hat{\gamma})} \theta (x,a_x) \Big)\hat{\lambda }( d \hat{\gamma}), \qquad \theta \in \Theta,
 \end{eqnarray}
 with $\hat{k}_\mu:\hat{\Gamma}_0 \to \mathbb{R}$ such that each
 $\hat{k}_\mu^{(n)}$, $n\in \mathbb{N}_0$, see (\ref{CS3}), has the
 following property: for Lebesgue-almost all $x_1, \dots, x_n\in X$,
$k_\mu^{(n)}$ defined by the expression
\begin{equation}
  \label{SF2}
k_\mu^{(n)} (x_1 , \dots , x_n) =
\int_{\mathbb{R}_{+}^n}\hat{k}_\mu^{(n)}((x_1,a_1), \dots,
(x_n,a_n))d a_1 \cdots d a_n.
\end{equation}
satisfies
\begin{equation}
  \label{SF2y}
  0\leq k_\mu^{(n)} (x_1 , \dots , x_n) \leq (n!)^\epsilon \varkappa^n.
\end{equation}
with certain $\epsilon \in [0,1)$ and $\varkappa>0$.
\end{definition}
For a given $\mu \in \mathcal{P}_{**}$, $\hat{k}_\mu$ and
$\hat{k}_\mu^{(n)}$ are called correlation function and $n$-th order
correlation function of $\mu$, respectively. It is worth noting that
$\hat{k}_\mu(\varnothing) =1$, which one readily gets from
(\ref{SF1}) with $F_\theta\equiv 1$, that corresponds to $\theta
\equiv 0$. Thus, by (\ref{SF2y}) each $k_\mu^{(n)}$ is a symmetric
element of $L^\infty(X^n)$. Note also that $\hat{k}^{(1)}_\mu(x,a)$
is the density of entities at point $x\in X$ and age $a\geq 0$. Then
$k^{(1)}_\mu(x)$ is merely the spatial density of entities. By
assuming that $k^{(1)}_\mu\in L^\infty (X)$ we allow the population
be infinite in state $\mu$, that holds if $k^{(1)}_\mu$ is not
integrable.

Let us show that $\mathcal{P}_{**}\subset \mathcal{P}_{*}$. By
standard formulas, for $\mu\in \mathcal{P}_{**}$, one gets, cf.
(\ref{J}), (\ref{SF2}) and (\ref{SF2y}),
\[
\mu(\Psi) = \int_{\hat{X}} \psi(x) \hat{k}_\mu^{(1)} (x, a) d x d a
= \int_{X} \psi(x) k_\mu^{(1)} (x, a) d x \leq \varkappa \int_{X}
\psi(x) d x.
\]
For $\mu\in \mathcal{P}_{**}$, its $q$-thinning amounts to
multiplying $\hat{k}_{\mu} (\hat{\eta})$ by $\prod_{x\in
p(\hat{\eta})} q(\hat{x})$. An important subclass of
$\mathcal{P}_{**}$ constitute Poisson measures
$\pi_{\hat{\varrho}}$. Each of them is completely determined by its
first-order correlation function $\hat{k}_\mu^{(1)}= \hat{\varrho}$
with $\hat{\varrho}(x,a)$ integrable in $a$ and essentially bounded
in $x$. In this case,
\begin{equation}
  \label{SF2u}
  \hat{k}_{\pi_{\hat{\varrho}}} (\hat{\eta}) = \prod_{x\in p(\hat{\eta})} \hat{\varrho}
  (\hat{x}),
\end{equation}
and hence (\ref{SF2y}) holds with $\epsilon=0$ and
\[
\varkappa = \esssup_{x\in X} \int_{\mathds{R}_{+}}
\hat{\varrho}(x,a) d a.
\]
Then by (\ref{SF1}) it follows that
\begin{equation}
  \label{CC}
  \pi_{\hat{\varrho}} (F_\theta) = \exp\left( \int_{\hat{X}} \hat{\varrho}(\hat{x}) \theta(\hat{x}) d \hat{x}
  \right).
\end{equation}
Note that each $\mu \in \mathcal{P}(\hat{\Gamma})$ can have a
correlation function understood as a distribution. To see this, let
us first define
\begin{equation*}
  %\label{D1}
 \delta (\hat{\xi}; \hat{\eta}) = \left\{\begin{array}{ll} \sum_{\sigma \in \varSigma_n}\prod_{j=1}^n \delta (\hat{x}_j -
 \hat{y}_{\sigma(j)}),
 \qquad &{\rm if} \  |\hat{\eta}|=|\hat{\xi}| = n; \\[.4cm] 0, \qquad &{\rm otherwise}.  \end{array} \right.
\end{equation*}
In the first line, $\hat{\xi}= \{\hat{x}_1, \dots , \hat{x}_n\}$,
$\hat{\eta}= \{\hat{y}_1, \dots , \hat{y}_n\}$, $\varSigma_n$ is the
symmetric group, and $\delta (\hat{x}-\hat{y})$ is the usual Dirac
$\delta$-function on $\mathds{R}^d \times \mathds{R}_{+}$. The
correlation function $\hat{k}_{\hat{\gamma}}$ of the
$\delta$-measure $\delta_{\hat{\gamma}}\in
\mathcal{P}(\hat{\Gamma})$ is then
\begin{equation}
  \label{D2}
  k_{\hat{\gamma}}(\hat{\xi}) = \sum_{\eta \subset p(\hat{\xi})} \delta (\hat{\xi};
  \hat{\eta}).
\end{equation}
By (\ref{SF1}) and (\ref{D2}) we then have
\[
\delta_{\hat{\gamma}} (F_\theta) = \int_{\hat{\Gamma}_0}
k_{\hat{\gamma}}(\hat{\xi}) \prod_{x\in p(\hat{\xi})} \theta
(\hat{x}) \hat{\lambda}( d \hat{\xi}) = \prod_{x\in p(\hat{\gamma})}
( 1 +\theta (\hat{x})).
\]
By means of $k_{\hat{\gamma}}$ one can define the correlation
function for any $\mu$ by the formula
\begin{equation}
  \label{D3}
  k_\mu (\hat{\xi}) = \int_{\hat{\Gamma}}
  k_{\hat{\gamma}}(\hat{\xi}) \mu (d \hat{\gamma}).
\end{equation}
Then for $\mu_1, \mu_2\in \mathcal{P} (\hat{\Gamma})$, by
(\ref{Conv1}) and (\ref{SF1}), and further by (\ref{Minlos}),
(\ref{D2}), (\ref{D3}), one readily gets that
\begin{equation*}
%\label{Conv2}
(\mu_1 \star \mu_2)(F_\theta) = \int_{\hat{\Gamma}_0}
\left(\sum_{\xi \subset p(\hat{\eta})} \hat{k}_{\mu_1}
(\hat{\eta}\setminus \hat{\xi}) \hat{k}_{\mu_2}(\hat{\xi}) \right)
\prod_{x\in p(\hat{\eta})} \theta (\hat{x})
\hat{\lambda}(d\hat{\eta}),
\end{equation*}
which by (\ref{SF2}) and (\ref{SF2y}) implies that $\mu_1 \star
\mu_2\in \mathcal{P}_{*}$, whenever $\mu_1 ,  \mu_2\in
\mathcal{P}_{*}$.

In the sequel, we use the Banach spaces
$\mathcal{G}_{\epsilon,\varkappa}$ with $\epsilon \in [0,1)$,
$\varkappa
>0$, of measurable functions $G:\hat{\Gamma}_0\to \mathbb{R}$ defined
by the following two properties. For each $n\in \mathds{N}$,
\begin{eqnarray}
  \label{SF3}
&(a)& \quad  |G|_n:= \esssup_{(x_1,x_2,\ldots,x_n)\in X^n}
\int_{(\mathbb{R}_+)^n} \left|G^{(n)}(x_1,a_1,\ldots, x_n,a_n)\right|da_1\ldots
da_n  < \infty \qquad \qquad \\[.2cm] \nonumber & (b)& \quad
\|G\|_{\epsilon,\varkappa}:= \sup_{n\in \mathbb{N}} |G|_n
(n!)^{-\epsilon} \varkappa^{-n} <\infty.
\end{eqnarray}
Note that $\hat{k}_\mu \in \mathcal{G}_{\epsilon,\varkappa}$ with
$\epsilon$ and $\varkappa$ as in (\ref{SF2y}). By (b) in (\ref{SF3})
one  concludes that
\begin{equation*}
%  \label{SF3a}
  \mathcal{G}_{\epsilon,\varkappa} \hookrightarrow \mathcal{G}_{\epsilon',\varkappa'}, \qquad
 \epsilon\leq \epsilon', \ \ \  \varkappa < \varkappa',
\end{equation*}
where $\hookrightarrow$ denotes continuous embedding.

\section{The Result}
In this section, we formulate and prove a statement describing the
evolution of our model. among others, we introduce the evolution
equations related to (\ref{KE}) and describe in which sense we are
going to solve them.
\subsection{The model and the result}

 The evolution of the considered population  is described by
(\ref{KE}) in which the Kolmogorov operator has the form
\begin{eqnarray}
  \label{S1}
 (LF)(\hat{\gamma}) & = & \sum_{x\in p(\hat{ \gamma})} \frac{\partial}{\partial a_x} F(\hat{\gamma})
  + \sum_{x\in p(\hat{\gamma})} m(\hat{x}) \left[F(\hat{\gamma} \setminus \hat{x}) - F(\hat{\gamma}) \right]
  \\[.2cm] \nonumber
 & + & \int_{X} b(x) \left[ F(\hat{\gamma} \cup (x,0)) - F(\hat{\gamma}) \right] d x, \nonumber
\end{eqnarray}
where the first term corresponds to aging, whereas the second and
third terms describe departing and arriving of the population
members, respectively. We assume that both $m$ and $b$ are
nonnegative, measurable and bounded. For further simplicity, with no
harm we additionally assume that $a\mapsto m(x,a)$ is continuous for
each $x\in X$. As mentioned above, we are not going to directly
solve the Kolmogorov equation. Instead, we consider the
corresponding Fokker-Planck equation\footnote{See \cite{FPER} for a
general theory of such equations}
\begin{equation}
  \label{FPE}
  \mu_t (F_\theta) = \mu_s(F_\theta) + \int_s^t \mu_u (L F_\theta)
  du, \qquad t> s \geq 0,
\end{equation}
for $F_\theta$ with $\theta\in \Theta$, see (\ref{SF}), where
$\Theta$ is the collection of all $\theta:\hat{X}\to \mathds{R}$
that have the following form
\begin{equation}
  \label{J10}
  \theta (x,a) = \vartheta (x) e^{-\tau \psi(x) \phi(a)} + e^{-\tau \psi(x)
  \phi(a)} -1.
\end{equation}
Here $\vartheta:X \to (-1, 0]$ is a continuous functions with
compact support, $\psi$ is as in (\ref{J}), $\tau\geq 0$ and
$\phi(a) = a/(1+a)$.  Let us then consider the collection
$\mathcal{F}= \{F_\theta: \theta \in \Theta\}$ with $F_\theta$
defined in (\ref{SF}), and hence of the form
\begin{equation*}
  %\label{J11}
F_\theta (\hat{\gamma}) = \exp\left(\sum_{x\in p(\hat{\gamma})} \log
(1+ \vartheta(x)) - \tau \sum_{x\in p(\hat{\gamma})} \psi(x)
\phi(a_x)\right).
\end{equation*}
Note that $0 < F_\theta (\hat{\gamma}) \leq 1$ for each
$\hat{\gamma}\in \hat{\Gamma}_*$ and $\mu(F_\theta) \leq 1$ for all
$\mu \in \mathcal{P}_*$. It is possible to show, cf. \cite[Theorem
18]{Dudley}, that each $F_\theta$ is $\upsilon_*$-continuous (see
(\ref{J2})). The pointwise product of $F_\theta$ and $F_{\theta'}$
is $F_{\theta''}$ with $\theta''$ corresponding to $\vartheta'' (x)=
\vartheta (x) + \vartheta'(x) + \vartheta(x) \vartheta'(x)$ and
$\tau'' = \tau + \tau'$. Assume that $\hat{\gamma}_1\neq
\hat{\gamma}_2$, both are in $\hat{\Gamma}_*$. Then one finds
$\hat{x}$ which belongs to exactly one of these configurations, say
$\hat{\gamma}_1$. If there is no $\hat{y}\in \hat{\gamma}_2$ with
$p(\hat{y})= p(\hat{x})$, one takes $\tau=0$ and $\vartheta$ such
that $\vartheta(p(\hat{x}))\neq 0$ and $\vartheta(p(\hat{y}))= 0$
for all $\hat{y}\in \hat{\gamma}_2$. Otherwise, one takes $\tau>0$
and $\vartheta(p(\hat{x})) = \vartheta(p(\hat{y})) \neq 0$ and
$\vartheta(p(\hat{z}))= 0$ for all $\hat{z}\in \hat{\gamma}_1\cup
\hat{\gamma}_2$ such that $\vartheta(p(\hat{z})) \neq
\vartheta(p(\hat{x}))$. In both cases, the corresponding $F_\theta$
separates $\hat{\gamma}_1$ and $\hat{\gamma}_2$, see property (c) of
Proposition \ref{J1pn}. Clearly, $F_\theta\equiv 1$ for $\tau=0$ and
$\vartheta \equiv 0$. Then by Proposition \ref{J1pn}
$\mathcal{F}=\{F_\theta: \theta \in \Theta\}$ is separating.

Let us prove now that $L F_\theta$ is $\mu$-integrable for each $\mu
\in \mathcal{P}_*$. By (\ref{S1}) we have
\begin{eqnarray}
  \label{Gua5}
  (L F_\theta)(\hat{\gamma}) & = & \sum_{x\in p(\hat{\gamma})} \left(
  \frac{\partial}{\partial a_x} \theta(x,a_x) - m(x,a_x)\theta(x,a_x)  \right) F_\theta(\hat{\gamma}\setminus
  \hat{x})\\[.2cm] \nonumber & + & F_\theta(\hat{\gamma}) \int_X b(x)
  \theta (x,0) d x =: H_1 (\hat{\gamma}) + H_2 (\hat{\gamma}).
\end{eqnarray}
Since $b$ is bounded, $H_2$ is also bounded. Since $\vartheta$ is
continuous and compactly supported, it is $\psi$-bounded. Hence, by
(\ref{J10}) one concludes that, for all $\hat{x}\in \hat{X}$, the
following holds
\begin{equation}
  \label{GS}
|\theta (\hat{x}) |\leq c_\theta \psi(x), \qquad
\left|\frac{\partial}{\partial a_x} \theta(\hat{x}) \right|\leq \tau
\psi (x),
\end{equation}
where $c_\theta$ depends only on the choice of $\vartheta$ and
$\tau$. By (\ref{GS}) we then have
\begin{equation}
  \label{GS1}
\left| H_1 (\hat{\gamma}) \right| \leq C_\theta \Psi (\hat{\gamma}),
\end{equation}
holding with an appropriate $C_\theta$. By Definition \ref{1adf}
this yields the property in question.

Now for $\theta \in \Theta$ and $m$ as in (\ref{SF}), we set
\begin{equation}
  \label{Gua}
  \theta_t (x,a) = \theta (x, a+t)\exp\left(- \int_a^{a+t} m(x, \alpha) d
  \alpha
  \right), \qquad t\geq 0,
\end{equation}
and then define a map $\mathcal{P}_*\ni \mu \mapsto \mu^t\in
\mathcal{P}_*$, $t\geq 0$ by the following relation
\begin{equation}
  \label{Gua1}
  \mu^t (F_\theta) = \mu (F_{\theta_t}), \qquad \theta \in \Theta.
\end{equation}
Since the family $\{F_\theta: \theta \in \Theta\}$ is separating,
each $\mu^t$ is uniquely determined by (\ref{Gua}), (\ref{Gua1}).
Note that the correlation function of $\mu^t$ can be expressed
through that of $\mu$ as follows
\begin{equation}
  \label{SFF}
 \hat{k}_{\mu^t} (\hat{\eta}) = \mathcal{J}_t (\hat{\eta}) \hat{k}_{\mu}
 (\hat{\eta}^t) \exp\left( - \sum_{x\in p(\hat{\eta})} \int^{a_x}_{a_x-t} m(x, \alpha) d
 \alpha\right),
\end{equation}
where $\hat{\eta}^t=\{ (x, a_x - t): x\in p(\hat{\eta}\}$,
\[
\mathcal{J}_t (\hat{\xi}) = \prod_{x\in p(\hat{\xi})} J_t(a_x),
\qquad J_t(a) := 1 - I_t(a),
\]
and $I_t(a):=\mathds{1}_{[0,t)}(a)$ is the indicator of $[0,t)$. By
(\ref{SFF}) the map $\mu \mapsto \mu^t$ preserves $\mathcal{P}_{**}$
and is a combination of a thinning and age shift. Now we are at a
position to formulate our result.
\begin{theorem}
  \label{1Tm}
For each $\mu_0 \in \mathcal{P}_*$,  the Fokker-Planck equation
(\ref{FPE}) has a solution of  the following form
\begin{equation}
  \label{Gua2}
  \mu_t = \mu_0^t \star \pi_{\hat{\varrho}_t},
\end{equation}
where $\mu_0^t$ is obtained from $\mu_0$ according to (\ref{Gua1})
and $\pi_{\hat{\varrho}_t}$ is the Poisson measure, see (\ref{SF2u})
and (\ref{CC}), with
\begin{equation}
  \label{Gua3}
  \hat{\varrho}_t (x,a) = \hat{b}(x,a) I_t(a) = b(x) \exp\left(
 - \int_0^a m(x, \alpha) d \alpha \right) I_t(a),
\end{equation}
and $I_t(a)$ being the indicator of $[0, t)$. If  $m(\hat{x}) \geq
m_*$ for some $m_*
>0$,  the evolution described in (\ref{Gua2}) has a stationary state
$\pi_{\hat{\varrho}}$ with $\hat{\varrho} = \hat{b}$, see
(\ref{Gua3}). In this case, the solution given in (\ref{Gua2}) with
$\mu_0 \in \mathcal{P}_{**}$ satisfies $\mu_t \Rightarrow
\pi_{\hat{\varrho}}$ as $t\to +\infty$, where we mean the usual weak
convergence of probability measures on $\hat{\Gamma}_*$.
\end{theorem}
Let us make some comments to this statement. According to
(\ref{Conv1}), (\ref{CC})  and (\ref{Gua2}) it follows that
\begin{equation}
  \label{Gua4}
  \mu_t (F_\theta) = \exp\left(\int_X \int_{[0,t)} \hat{b}(x,a) \theta (x,a) d x d a
  \right)\mu_0 (F_{\theta_t}).
\end{equation}
Hence, the solution satisfies the initial condition $\mu_t|_{t=0} =
\mu_0$, see (\ref{Gua}).  If $\mu_0 (\varnothing) = 1$, i.e., the
initial state is an empty habitat, by (\ref{Gua4}) it follows that
$\mu_t = \pi_{\hat{\varrho}_t}$ with $\hat{\varrho}_t$ given in
(\ref{Gua3}). Let us show that this $\mu_t$ satisfies (\ref{FPE}).
For a Poisson measure $\pi_{\hat{\varrho}}$, by (\ref{Minlos}) and
(\ref{SF}) we have that
\begin{gather}
\label{Gua6} \pi_{\hat{\varrho}} (H_1) = \int_{\hat{\Gamma}_0}
\left(\prod_{x\in p(\hat{\eta})}\hat{ \varrho} (\hat{x} )
\right)\sum_{x\in p(\hat{\eta})} \left[ \frac{\partial}{\partial
a_x} - m(x,a_x)\right] \theta(\hat{x}) \prod_{y \in
p(\hat{\eta}\setminus \hat{x})} \theta(\hat{y}) \hat{\lambda}( d
\hat{\eta})\\[.2cm] \nonumber = \int_{\hat{\Gamma}_0}\left(\prod_{x\in
p(\hat{\eta})} \hat{\varrho}(\hat{x}) \right)\left( \int_{\hat{X}}
\hat{\varrho}(\hat{x})\left[ \frac{\partial}{\partial a_x} -
m(x,a_x)\right] \theta(\hat{x}) d \hat{x} \right)\prod_{y \in
p(\hat{\eta})} \theta(\hat{y}) \hat{\lambda}( d \hat{\eta})\\[.2cm] \nonumber
= - \left(\int_X \hat{\varrho}(x,0 ) \theta(x,0) d  x \right)
\pi_{\hat{\varrho}}(F_\theta) \\[.2cm] \nonumber - \left( \int_{\hat{X}}
\theta(\hat{x})\left[ \frac{\partial}{\partial a_x} +
m(x,a_x)\right] \hat{\varrho}(\hat{x}) d \hat{x} \right)
\pi_{\hat{\varrho}}(F_\theta).
\end{gather}
And also
\begin{equation}
  \label{Gua7}
\pi_{\hat{\varrho}} (H_2) = \left(\int_X  b(x) \theta(x,0 ) d x
\right) \pi_{\hat{\varrho}}(F_\theta).
\end{equation}
In the sense of distributions, we have that
\begin{equation*}
  %\label{Gua7a}
\frac{\partial}{\partial a} I_t (a) = - \frac{\partial}{\partial t}
I_t(a).
\end{equation*}
Then for $\hat{\varrho}_t(\hat{x})$ given in (\ref{Gua3}), one
obtains
\begin{equation}
  \label{Gua7b}
\left[ \frac{\partial}{\partial a_x} + m(x,a_x)\right]
\hat{\varrho}_t(\hat{x}) = - \frac{\partial}{\partial
t}\hat{\varrho}_t(\hat{x}).
\end{equation}
By (\ref{Gua6}), (\ref{Gua7}) and the latter equality it follows
that
\begin{eqnarray}
  \label{Gua8}
\pi_{\hat{\varrho}_t} (L F_\theta) & = &  \exp\left(\int_{\hat{X}}
\hat{\varrho}_t(\hat{x}) \theta (\hat{x})  d\hat{x} \right)
\frac{\partial}{\partial t} \int_{\hat{X}} \hat{\varrho}_t(\hat{x})
\theta (\hat{x})  d\hat{x} \\[.2cm] \nonumber & = & \frac{\partial}{\partial t}\exp\left(\int_{\hat{X}}
\hat{\varrho}_t(\hat{x}) \theta (\hat{x})  d\hat{x} \right),
\end{eqnarray}
by which one readily concludes that $\mu_t = \pi_{\hat{\varrho}_t}$
satisfies (\ref{FPE}).

\subsection{Proof of Theorem \ref{1Tm}}
The proof of the first part will be done by showing that: (a) for
each $\theta\in \Theta$, the map $t \mapsto \mu_t (F_\theta)$ has a
continuous derivative at each $t>0$; (b) this derivative satisfies,
cf. (\ref{Gua8}),
\begin{equation}
  \label{GS2}
  \frac{d}{dt} \mu_t (F_\theta) = \mu_t(LF_\theta).
\end{equation}
By (\ref{Conv1}), (\ref{Gua2}), (\ref{Gua3}) and (\ref{Gua4}) we
have
\begin{eqnarray}
  \label{GS3}
\mu_t (F_\theta)& = & \mu_0 (F_{\theta_t}) \pi_{\hat{\varrho}_t}
(F_\theta)=:   \mu_0 (F_{\theta_t}) Q_\theta (t).
\end{eqnarray}
In view of (\ref{Gua8}), the continuous differentiability in
question will thus follow by the same property of $t\mapsto \mu_0
(F_{\theta_t})$. By (\ref{Gua}) we have
\begin{eqnarray}
  \label{GS4}
  \frac{\partial }{\partial t} F_{\theta_t} (\hat{\gamma}) & = &
\sum_{x\in p(\hat{\gamma})} \left( \frac{\partial }{\partial a_x}
\theta_t(\hat{x})\right) F_{\theta_t} (\hat{\gamma}\setminus
\hat{x}) \\[.2cm] \nonumber & - & \sum_{x\in p(\hat{\gamma})} m(x, a_x ) \theta_t
(\hat{x}) F_{\theta_t} (\hat{\gamma}\setminus \hat{x})\\[.2cm] \nonumber & =:&
\sum_{x\in p(\hat{\gamma})} \sigma_t (\hat{x} )F_{\theta_t}
(\hat{\gamma}\setminus \hat{x}) =: S_t (\hat{\gamma}).
\end{eqnarray}
Similarly as in (\ref{GS1}) we then conclude that
\[
\left|   \frac{\partial }{\partial t} F_{\theta_t} (\hat{\gamma})
\right| \leq C'_\theta \Psi(\hat{\gamma}),
\]
with a certain $C'_\theta>0$. By Lebesgue's dominated convergence
theorem and (\ref{SF1}) this yields
\begin{eqnarray}
  \label{GS5}
 \frac{d}{dt} \mu_0 (F_{\theta_t})& = & \mu_0 \left( \frac{\partial }{\partial t}
 F_{\theta_t}\right)=\mu_0 (S_t)  \\[.2cm] \nonumber & = & \int_{\hat{\Gamma}_0} \hat{k}_{\mu_0} (\hat{\eta}) \sum_{x\in p(\hat{\eta})}\left( \sigma_t(\hat{x})
\prod_{y\in p(\hat{\eta}\setminus \hat{x})} \theta_t (\hat{y})
\right) \hat{\lambda}(d \hat{\eta}) ,
\end{eqnarray}
as well as the continuity of the map $t \mapsto \mu_0 \left(
\frac{\partial }{\partial t} F_{\theta_t}\right)$. Here
$\hat{k}_{\mu_0}$ is the correlation function of $\mu_0$ understood
in the sense of (\ref{D3}). Now let us turn to proving (\ref{GS2}).
By (\ref{GS3}) and (\ref{GS5}) we have
\begin{eqnarray}
  \label{GS6}
  {\rm LHS}(\ref{GS2}) =  \mu_0 (S_t) Q_\theta (t) +  \mu_t (F_\theta) \int_X b(x) \theta (x,t) \exp\left( - \int_0^t m(x,\alpha) d
  \alpha\right) d x.
\end{eqnarray}
At the same time, by (\ref{Gua5}) it follows that
\[
H_1 (\hat{\gamma}_1 \cup \hat{\gamma}_2) = H_1 (\hat{\gamma}_1 )
F_\theta (\hat{\gamma}_2) + H_1 (\hat{\gamma}_2 ) F_\theta
(\hat{\gamma}_1),
\]
which by (\ref{Conv}) and (\ref{Gua2}) yields
\begin{eqnarray}
  \label{GS7}
  {\rm RHS}(\ref{GS2}) = \mu_0^t (H_1) \pi_{\hat{\varrho}_t}
  (F_\theta) + \mu_0^t (F_\theta)  \pi_{\hat{\varrho}_t}(H_1)
+ \mu_t (F_\theta) \int_X b(x) \vartheta (x) d x, \qquad
\end{eqnarray}
Note that
\begin{equation}
  \label{GS7a}
\pi_{\hat{\varrho}_t}
  (F_\theta) = \exp\left( \int_{\hat{X}} \hat{\varrho}_t (\hat{x} )\theta(\hat{x}) d \hat{x}
  \right) = Q_\theta (t),
\end{equation}
see (\ref{CC}), (\ref{Gua3}) and (\ref{GS3}). By (\ref{Gua5}) we
have that
\begin{eqnarray*}
%  \label{GS8}
  H_1 (\hat{\gamma} ) &  = & \sum_{x\in p(\hat{\gamma})} h_\theta
  (\hat{x}) F_\theta (\hat{\gamma}\setminus \hat{x}),
  \\[.2cm]\nonumber h_\theta (x,a) & := & \frac{\partial}{\partial a}
  \theta (x,a) - m(x,a) \theta(x,a).
\end{eqnarray*}
By (\ref{SF1}), (\ref{Minlos})  and (\ref{SFF}) one then gets
\begin{eqnarray}
  \label{GS9}
  \mu_0^t (H_1) & = & \int_{\hat{\Gamma}_0} \hat{k}_{\mu_0^t}
  (\hat{\eta}) \left( \sum_{x\in p(\hat{\eta})} h_\theta (\hat{x}) \prod_{y\in p(\hat{\eta}\setminus \hat{x})} \theta (\hat{y})
  \right) \hat{\lambda}(d \hat{\eta}) \\[.2cm] \nonumber & =
  &\int_{\hat{\Gamma}_0} \left( \int_{\hat{X}} \hat{k}_{\mu_0^t}
  (\hat{\eta}\cup \hat{x}) h_\theta (\hat{x} )d \hat{x} \right) \prod_{y\in p(\hat{\eta})} \theta (\hat{y}) \hat{\lambda}(d \hat{\eta})
  \\[.2cm] \nonumber & =
  & \int_{\hat{\Gamma}_0} \left( \int_{X} \int_t^{+\infty}\hat{k}_{\mu_0} (\hat{\eta}^t \cup (x,a-t)) \exp\left( -\int_{a-t}^{a}
  m(x,\alpha) d \alpha \right) h_\theta (x,a) d x d a \right) \\[.2cm] \nonumber
  & \times & \mathcal{J}_t (\hat{\eta}) \left( \prod_{y\in p(\hat{\eta})}
  \theta (y,a_y) \exp\left( - \int_{a_y-t}^{a_y} m(y , \alpha) d \alpha
  \right)\right)\hat{\lambda} ( d \hat{\eta})  \\[.2cm] \nonumber
  & = &  \int_{\hat{\Gamma}_0}  K_t (\hat{\eta})  \prod_{y\in p(\hat{\eta})}
  \theta_t (\hat{y} ) \hat{\lambda} ( d \hat{\eta}).
\end{eqnarray}
Here $\hat{\eta}^t$ and $\mathcal{J}_t$ are as in  (\ref{SFF}) and
$\theta_t$ is defined in (\ref{Gua}), whereas
\begin{gather*}
 K_t (\hat{\eta}) := \int_{X} \int_{t}^{+\infty}
 \hat{k}_{\mu_0}(\hat {\eta}^t \cup(x,a-t)) h_\theta (x,a) \exp\left(
 -\int_{a-t}^a m(x, \alpha) d \alpha\right) d x d a \\[.2cm]
 \nonumber = \int_{X} \int_{0}^{+\infty} \hat{k}_{\mu_0}(\hat {\eta} \cup(x,a)) h_\theta (x,a+t) \exp\left(
 -\int_{a}^{a+t} m(x, \alpha) d \alpha\right) d x d a.
\end{gather*}
By (\ref{Gua}) and (\ref{GS4}) we have
\[
h_\theta (x,a+t) \exp\left(
 -\int_{a}^{a+t} m(x, \alpha) d \alpha \right)= \sigma_t (x,a).
\]
We use this in the latter expression and then in (\ref{GS9}) and
thus arrive at the following
\begin{equation}
  \label{GS10}
  \mu_0^t (H_1) = \int_{\hat{\Gamma}_0} \hat{k}_{\mu_0} (\hat{\eta})
 \left( \sum_{x\in p(\hat{\eta})} \sigma_t(\hat{x}) \prod_{y \in
 p(\hat{\eta}\setminus \hat{x})} \theta_t (\hat{y}) \right)
 \hat{\lambda}(d \hat{\eta}) = \mu_0(S_t),
\end{equation}
see (\ref{GS5}). Now similarly as in (\ref{Gua6}) we obtain
\begin{eqnarray}
  \label{GS11}
  \pi_{\hat{\varrho}_t}(H_1)& = & \left(\int_{\hat{X}} \hat{\varrho}_t(\hat{x}) h_\theta (\hat{x}) d \hat{x}
  \right)\pi_{\hat{\varrho}_t}(F_\theta) \\[.2cm] \nonumber & = &  \left(
  \int_X \int_0^t b(x) \exp\left( -\int_0^a m(x,\alpha) d
  \alpha\right) \left[ \frac{\partial}{\partial a} - m(x,a)\right] \theta(x,a) d x d a
  \right)\pi_{\hat{\varrho}_t}(F_\theta) \\[.2cm] \nonumber & = &  \left(
\int_X b(x) \left[ \exp\left( -\int_0^t m(x,\alpha) d \alpha
\right)\theta ( x , t) - \theta (x,0)\right] d x \right)
\pi_{\hat{\varrho}_t}(F_\theta) .
\end{eqnarray}
Finally, we use (\ref{GS10}) and (\ref{GS11}) in (\ref{GS7}), take
into account (\ref{GS7a}) and (\ref{GS6}),  and conclude that
(\ref{GS2}) holds true.

To prove that $\pi_{\hat{\varrho}}$ with $\hat{\varrho}=\hat{b}$ is
a stationary solution of (\ref{FPE}) we again use (\ref{Gua5}) and
(\ref{Gua6}). For
\[
\hat{\varrho}(x,a)=\hat{b}(x,a) = b(x) \exp\left( - \int_0^a m(x,
\vartheta) d \vartheta\right),
\]
we have, cf. (\ref{Gua7b}),
\[
\left[\frac{\partial}{\partial a_x} + m(x,a_x) \right]
\hat{b}(x,a_x) =0,
\]
which by (\ref{Gua6}) yields $\pi_{\hat{\varrho}}(L F_\theta)=0$,
and hence the property in question.

To complete the proof we have to show the stated weak convergence,
assuming $m(x,a) \geq m_* >0$ and $\mu_0\in \mathcal{P}_*$. The
latter fact implies $\hat{k}_{\mu_0}\in
\mathcal{G}_{\epsilon,\varkappa}$ for some $\epsilon \in [0,1)$ and
$\varkappa >0$, see Definition \ref{1df}. Recall that we also assume
that $b(x)\leq b^*$ for some $b^*>0$. Let us prove that $\hat{k}_t
\in \mathcal{G}_{\epsilon, \varkappa^*}$ with the same $\epsilon$
and $\varkappa^* = \max\{\varkappa; b^*/m_* \}$. Then for each $l\in
\mathds{N}_0$, we have that
\begin{equation}
  \label{TT}
 0\leq  k_{\mu_0}(x_1, \dots , x_l) \leq (l!)^\epsilon \varkappa^l,
\end{equation}
see (\ref{SF2}), (\ref{SF2y}). For a given $n\in \mathds{N}$, by
(\ref{TT}) we obtain
\begin{eqnarray}
  \label{TT1}
 k_t^{(n)} (x_1, \dots , x_n) & \leq & \sum_{l=0}^n {n \choose l}\left(
 b^* \int_0^t e^{- m_* a} d a \right)^{n-l}e^{- l
 m_*t}\\[.2cm] \nonumber & \times &\int_{[t,+\infty)^l} \hat{k}_{\mu_0} ((x_1, a_1-t), \dots ,
 (x_l, a_l-t)) d a_1 \cdots  d a_l \\[.2cm] \nonumber &\leq & (n
 !)^\epsilon \sum_{l=0}^n {n \choose l} \left( \frac{b^*}{m_*} [1 - e^{- m_* t}]\right)^{n-l}
 \left( \varkappa e^{- m_* t} \right)^l \\[.2cm] \nonumber & \leq & (n
 !)^{\epsilon} (\varkappa^*)^n.
\end{eqnarray}
Thus, $\hat{k}_t \in \mathcal{G}_{\epsilon, \varkappa^*}$ for all
$t\geq 0$. Let us prove that
\begin{equation}
  \label{TT2}
  \| \hat{k}_t  - \hat{k}_{\pi_{\hat{\varrho}}}\|_{\epsilon,\kappa}
  \to 0, \qquad {\rm as}  \ \ t\to +\infty,
\end{equation}
for some $\kappa\geq \varkappa^{*}$. To this end, for some $t>0$, we
write
\begin{eqnarray*}
  %\label{TT3}
\hat{k}_{\pi_{\hat{\varrho}}} (\hat{\eta}) & = & \prod_{y\in
p(\hat{\eta})} b(y) \exp\left( - \int_0^{a_y} m(y,\vartheta) d
\vartheta \right)(I_t (a_y) + J_t(a_y))\\[.2cm] \nonumber & = &
\sum_{\xi\subset p(\hat{\eta})} \hat{k}_{\pi_{\hat{\varrho}}}
(\hat{\xi})\mathcal{J}_t(\hat{\xi}) \prod_{x\in
p(\hat{\eta}\setminus \hat{\xi})} \hat{b}(\hat{x}) I_t (a_x)
\end{eqnarray*}
Next -- similarly as in (\ref{TT1}) -- for a given $n\in
\mathds{N}$, we have, see (\ref{SF3}),
\begin{eqnarray}
  \label{TT4}
\left|\hat{k}_t  - \hat{k}_{\pi_{\hat{\varrho}}} \right|_n & \leq &
\sum_{l=0}^{n-1}{n \choose l} (b^*)^l \left( \int_0^t e^{- m_* a} d
a \right)^l e^{-(n-l) m_* t} \\[.2cm] \nonumber & \times & \left[
\left(\frac{b^*}{m_*} \right)^{n-l} + ((n-l)!)^\epsilon
\varkappa^{n-l} \right] \\[.2cm] \nonumber & \leq & 2  (n!)^\epsilon
(\varkappa^*)^n\sum_{l=0}^{n-1}{n \choose l}  e^{- (n-l)m_* t} \\[.2cm] \nonumber & \leq &
(n!)^\epsilon n ( 2 \varkappa^*)^n e^{- m_* t} \leq (n!)^\epsilon (
3 \varkappa^*)^n e^{- m_* t}.
\end{eqnarray}
By the latter estimate and (\ref{SF3}) one gets (\ref{TT2}) holing
with $\kappa = 3 \varkappa^{*}$. Note that (\ref{TT2}) implies that
\begin{equation}
  \label{TT2a}
 \forall \theta \in \Theta \qquad \mu_t (F_\theta) \to
 \pi_{\hat{\varrho}}(F_\theta) = \exp\left(\int_{\hat{X} } \hat{b}(\hat{x}) \theta ( \hat{x}) d \hat{x}
 \right), \quad t \to +\infty.
\end{equation}
Indeed, let $\theta \in \Theta$ and a compact $\varLambda$ be such
that $\theta (x,a) =0$ for $x \in \varLambda^c$.  By (\ref{SF1}) and
(\ref{TT4}) we then have
\begin{gather*}
\left| \mu_t (F_\theta) -
 \pi_{\hat{\varrho}}(F_\theta)\right| \leq \int_{\hat{\Gamma}_0}
 \left| \hat{k}_t (\hat{\eta}) -
 \hat{k}_{\pi_{\hat{\varrho}}}(\hat{\eta})\right| \left( \prod_{x\in p(\hat{\eta})} \left|\theta (x, a_x)\right|\right)
 \hat{\lambda}(d \hat{\eta})\\[.2cm] \leq e^{- m_8 t} \sum_{n=1}^\infty
 \frac{\left( \kappa {\rm
 Vol}(\varLambda)\right)^n}{(n!)^{1-\epsilon}},
\end{gather*}
which yields (\ref{TT2a}).

Now we show that the family $(\mu_t)_{t\geq 0}$ is tight, which by
Prohorov's theorem would yield its relative weak compactness, and
hence the existence of of the corresponding accumulation points.
Similarly as in \cite[Corollary A2.6V, page 406]{DV1} one proves
that a subset, $\hat{\Delta} \subset \hat{\Gamma}$, is compact (in
the vague topology) if, for each compact $\varLambda \subset X$,
there exists $c_\varLambda >0$ such that
\begin{equation*}
%\label{TT10}
\forall \hat{\gamma}\in \hat{\Delta} \qquad
N_\varLambda (\hat{\gamma}):=|p(\hat{\gamma}) \cap \varLambda|\leq
c_{\varLambda}.
\end{equation*}
Then the tightness in question will follow by the fact that
\begin{equation}
  \label{TT11}
\forall t\geq 0 \qquad \mu_t(N_\varLambda) \leq C_\varLambda,
\end{equation}
holding for some $C_\varLambda >0$. At the same time, for $\mu\in
\mathcal{P}_*$, we know that
\begin{equation*}
  %\label{TT12}
  \mu (N_\varLambda) = \int_{\hat{X}} k_\mu^{(1)} (\hat{x})
  \mathds{1}_\varLambda (x) d \hat{x} \leq {\rm Vol}(\varLambda)
  |\hat{k}_\mu|_1,
\end{equation*}
where ${\rm Vol}(\varLambda)$ is Lebesgue's measure (volume) of
$\varLambda$ and $|\hat{k}_\mu|_1$ is defined in (\ref{SF3}). By
(\ref{TT1}) we then obtain $\mu_t (N_\varLambda) \leq {\rm
Vol}(\varLambda)\varkappa^*$, which yields (\ref{TT11}) and thereby
the relative weak compactness of $(\mu_t)_{t\geq 0}$. Let $\mu\in
\mathcal{P}_*$ be an accumulation point of $(\mu_t)_{t\geq 0}$, and
hence $\mu_{t_n}\Rightarrow \mu$ for some $\{t_n\}_{n\in
\mathds{N}}$, $t_n \to +\infty$. Therefore,
\[
\mu_{t_n} (F_\theta) \to \mu(F_\theta), \qquad n\to +\infty,
\]
for each $\theta \in \Theta$. By (\ref{TT2a}) this immediately
yields $\mu= \pi_{\hat{\varrho}}$ as the family $\{ F_\theta: \theta
\in \Theta\}$ is separating. This completes the proof of the whole
theorem.

\section*{Acknowledgements}
The research of both authors was financially supported by National
Science Centre, Poland, grant 2017/25/B/ST1/00051, that is cordially
acknowledged by them.

\end{document}